\magnification=1200
\def\qed{\unskip\kern 6pt\penalty 500\raise -2pt\hbox
{\vrule\vbox to 10pt{\hrule width 4pt\vfill\hrule}\vrule}}
\centerline
	{GRAPH-COUNTING POLYNOMIALS FOR ORIENTED GRAPHS.}
\bigskip
\centerline{by David Ruelle\footnote{$\dagger$}{Math. Dept., Rutgers University, and 
IHES, 91440 Bures sur Yvette, France. email: ruelle@ihes.fr}.}
\bigskip\bigskip\bigskip\bigskip\noindent
	{\leftskip=1.8cm\rightskip=2cm{\sl Abstract:}
\medskip\noindent
If ${\cal F}$ is a set of subgraphs $F$ of a finite graph $E$ we define a graph-counting polynomial
$$	p_{\cal F}(z)=\sum_{F\in{\cal F}}z^{|F|}   $$
In the present note we consider oriented graphs and discuss some cases where ${\cal F}$ consists of unbranched subgraphs $E$.  We find several situations where something can be said about the location of the zeros of $p_{\cal F}$.\par}
\vfill\eject
	Let ${\cal F}$ be a set of subgraphs $F$ of a finite graph $E$.  We denote by $|F|$ the number of edges of $F$ and define a polynomial
$$	p_{\cal F}(z)=\sum_{F\in{\cal F}}z^{|F|}   $$
(graph-counting polynomial associated with ${\cal F}$).  The case of unoriented graphs has been discussed earlier (see [4]-[6] and [1]-[3]); here we mostly consider oriented graphs.
\medskip
	We shall find that for suitable ${\cal F}$ we can restrict the location of the zeros of $p_{\cal F}$ (for instance to the imaginary axis).  The proofs will be based on the following fact:
\medskip
	{\it Lemma} (Asano-Ruelle).  {\sl Let $K_1,K_2$ be closed subsets of the complex plane ${\bf C}$ such that $K_1,K_2\not\ni0$ and assume that
$$	A+Bz_1+Cz_2+Dz_1z_2\ne0\qquad{\rm when}\qquad z_1\notin K_1,z_2\notin K_2   $$
Then
$$	A+Dz\ne0\qquad{\rm when}\qquad z\notin-K_1K_2   $$
where $-K_1K_2$ is minus the set of products of an element of $K_1$ and an element of $K_2$.}  (The replacement of $A+Bz_1+Cz_2+Dz_1z_2$ by $A+Dz$ is called Asano contraction).
\medskip
	For a proof see for instance the Appendix A of [6].  The results given below follow rather directly from this lemma.
\bigskip\noindent
{\bf 1. Definitions.  Subgraphs of an oriented graph.}
\medskip
	We say that a pair $(V,E)$ of finite sets is an {\it oriented graph} if $V\ne\emptyset$ and we are given two maps $x',x'':E\to V$ such that $x'(e)\ne x''(e)$ for all $e\in E$.  The elements $x$ of $V$ are {\it vertices}, and the elements $e$ of $E$ are {\it oriented edges} with endpoints $x'(e),x''(e)$; $e$ is outgoing at $x'(e)$ and ingoing at $x''(e)$ and we write $e: x'(e)\to x''(e)$.  We allow different edges $e_1,e_2$ such that $e_1:x'\to x''$ and $e_2:x'\to x''$ or $e_2:x''\to x'$.
\medskip
	We say that $(V,E)$ is {\it bipartite} if we are given a partition of $V$ into nonempty sets $V_1$  and $V_2$ such that for each $e\in E$ the points $x'(e)$ and $x''(e)$ are in different sets of the partition $\{V_1,V_2\}$.  If $e\in E$ we have thus either $e:V_1\to V_2$ or $e:V_2\to V_1$.
\medskip
	We call a subset $F$ of $E$ a {\it subgraph} of $(V,E)$.  We say that $F$ is connected if for each partition $\{F_1,F_2\}$ of $F$ there is an $x\in V$ which is an endpoint of both some $e_1\in F_1$ and some $e_2\in F_2$.  A subgraph is thus a union of connected components $F_j$ in a unique way.	
\medskip
	We say that $F$ is an {\it unbranched} subgraph if, for each $x\in V$,
$$	|\{e\in F:x'(e)=x\}|\le1\qquad{\rm and}\qquad|\{e\in F:x''(e)=x\}|\le1   $$
We say that an unbranched subgraph is a {\it loop subgraph} if for each $x\in V$ we have
$$	|\{e\in F:x'(e)=x\}|=|\{e\in F:x''(e)=x\}|   $$
We denote by $U(E)$, resp. $L(E)$, the set of unbranched subgraphs, resp. loop subgraphs of the oriented graph $(V,E)$.  If $F$ is an unbranched subgraph, we can write $F$ as a disjoint union of connected components $F_j$ which are either loops (i.e., $F_i\in L(E)$) or, if they are not loops, have different endpoints $x'_j$ and $x''_j$ such that $x'_j\to\cdots(e)\cdots\to x''_j$ for each $e\in F_j$.
\medskip
	Introduce now complex variables $z'_e,z''_e$, write $Z'=(z'_e)_{e\in E}$, $Z''=(z''_e)_{e\in E}$ and, for each $x\in V$, let
$$	p_x(Z',Z'')=\big(a'(x)+\sum_{e:x'(e)=x}z'_e\big)\big(a''(x)+\sum_{e:x''(e)=x}z''_e\big)\eqno{(1)}   $$
with some choice of the $a'(x),a''(x)\in{\bf C}$.  For small $\epsilon>0$ we also write
$$	p_x^\epsilon(Z',Z'')
	=\big(a'(x)+\sum_{e:x'(e)=x}(z'_e+\epsilon)\big)(a''(x)+\sum_{e:x''(e)=x}\big(z''_e+\epsilon)\big)   $$
and
$$	\tilde p_x=1+p_x\qquad,\qquad\tilde p_x^\epsilon=1+p_x^\epsilon   $$
Choosing between $p_x$ and $\tilde p_x$ for each $x$ and applying Asano contractions $(z'_e,z''_e)\to z_e$ for all $e\in E$ to the polynomials
$$	\prod_{x\in V}(p_x(Z',Z'')\hbox{ or }\tilde p_x(Z',Z''))\qquad,\qquad
	\prod_{x\in V}(p_x^\epsilon(Z',Z'')\hbox{ or }\tilde p_x^\epsilon(Z',Z''))\eqno{(2)}   $$
we obtain polynomials
$$	P(Z)\qquad,\qquad P^\epsilon(Z)   $$
where $Z=(z_e)_{e\in E}$ and
$$	\lim_{\epsilon\to0}P^\epsilon(Z)=P(Z)   $$
We shall obtain examples of $p(z)=p_{\cal F}(z)$ by taking all components $z_e$ of $Z$ equal to $z$.
\bigskip\noindent
{\bf 2. Unbranched subgraphs of an oriented graph.}
\medskip
If there are only $p_x$ factors in $(2)$ and we assume $a'(x)a''(x)=1$ for all $x$, we have
$$	P(Z)=\sum_{F\in U(E)}\prod_{F_j\subset F}^{\rm conn}
	[a''(x'_j)a'(x''_j)]_{F_j{\rm not}\,{\rm loop}}\prod_{e\in F_j}z_e\eqno{(3)}   $$
where the product is over the connected components $F_j$ of $F$ and $x'_j$, $x''_j$ are the endpoints of $F_j$ if $j$ is not a loop; $[a''(x'_j)a'(x''_j)]$ is replaced by 1 if $F_j$ is a loop.
\medskip
	If we take $a'(x)=a''(x)=1$ and set all $z_e$ equal to $z$ we obtain the {\it unbranched subgraph counting polynomial}
$$	p_{\rm unbranched}(z)=\sum_{F\in U(E)}z^{|F|}\eqno{(4)}   $$
\indent
	{\bf 2.1. Proposition.}  {\sl The zeros of the unbranched subgraph counting polynomial (4) are all real $<0$.}
\medskip
	To prove this let $\alpha', \alpha''\in(-\pi/2,\pi/2)$.  Assuming
$$	{\rm Re}(z'_e+\epsilon)e^{-i\alpha'}>0\qquad,\qquad{\rm Re}(z''_e+\epsilon)e^{-i\alpha''}>0   $$
for all $e\in E$, we have $\prod_{x\in V}p_x^\epsilon(Z',Z'')\ne0$ and therefore by Asano contraction $P^\epsilon(Z)\ne0$ if $e^{-i(\alpha'+\alpha'')}z_e$ is in a neighborhood of the positive real axis and $-\pi<\alpha'+\alpha''<\pi$.  Let $p(z)$ and $p^\epsilon(z)$ be obtained by taking all $z_e$ equal to $z$ in $P(Z)$ and $P^\epsilon(Z)$.  Then $p^\epsilon(z)\ne0$ if $\arg z\in(-\pi,\pi)$.  Using Hurwitz's theorem we let $\epsilon\to0$ in $p^\epsilon(z)$ and find that either $p(z)$ vanishes identically or $p(z)\ne0$ if $\arg z\in(-\pi,\pi)$.  Clearly $p(0)\ne0$ because $\emptyset\subset U(E)$ and we obtain thus $p_{\rm unbranched}(z)=p(z)\ne0$ if $z$ is not real $<0$ .\qed
\medskip\noindent
[In fact since $a'(x)=a''(x)=1$ we could have done without $\epsilon$ in the present situation].
\medskip
	Let now deg$_2$ be the max over $x\in V$ of the number $d'_x$ of outgoing edges at $x$ times the max over $x$ of the number $d''_x$ of ingoing edges at $x$.  Then $p_x(Z',Z'')\ne0$ if $|z'_e|<1/d'_x$ and $|z''_e|<1/d''_x$ for all $e\in E$, so that $P(Z)\ne0$ if all $z_e<1/{\rm deg}_2$.  Finally $p(z)\ne0$ if $|z|<1/{\rm deg}_2$, i.e., the zeros of $p(z)$ are negative $\le-1/{\rm deg}_2$.
\medskip
	{\bf 2.2. Remark.}
\medskip
	Given $V_0\subset V$ let $a'(x)=a''(x)=1$ if $x\in V_0$ and $a'(x)=a''(x)=0$ if $x\notin V_0$.  The polynomial $p$ counts then unbranched polynomials going through all $x\notin V_0$ and one finds that the zeros of $p$ are real $\le0$  if $p$ does not vanish identically.
\medskip
	The following result is relevant to Section 3.2 below.
\medskip
	{\bf 2.3. Proposition.}  {\sl Let $(V,E)$ be bipartite corresponding to a partition $\{V_1,V_2\}$ of $V$ and let $U_{\rm even}(E)$ consist of the unbranched subgraphs $F$ such that $|F_j|$ is even for each connected component $F_j$ of $F$.  We define
$$	p_{\rm unbranched\,even}(z)=\sum_{F\in U_{\rm even}(E)}z^{|F|}\eqno{(5)}   $$
Assume that the connected subgraphs in $U(E)$ come in pairs $(G,\bar G)$ connecting the same vertices and both $G$, $\bar G$ are loops or not loops.  If $x',x''$ (resp. $\bar x',\bar x''$) are the endpoints of non-loop $G$ (resp. $\bar G$) we also assume $x'=\bar x''$, $x''=\bar x'$.  Under these conditions the zeros of the polynomial (5) are all purely imaginary.}
\medskip
	In equation (1) we take $a'(x)=(1+i)/\sqrt2$ if $x\in V_1$, $a'(x)=(1-i)/\sqrt2$ if $x\in V_2$ and let $a''(x)$ be the complex conjugate $a'(x)^*$ of $a'(x)$ in all cases.  We use thus the $p_x(Z',Z'')$, $p_x^\epsilon(Z',Z'')$ corresponding to those $a'(x)$, $a''(x)$.
\medskip
	We obtain polynomials $P(Z)$, resp. $P^\epsilon(Z)$, by Asano contractions of $\prod p_x(Z',Z'')$, resp. $\prod p_x^\epsilon(Z',Z'')$, and (3) gives
$$	P(Z)=\sum_{F\in U(E)}\prod_{F_j\subset F}^{\rm conn}\gamma(F_j)\prod_{e\in F_j}z_e\eqno{(6)}   $$
with $\gamma(F_j)=a''(x'_j)a'(x''_j)$ if $F_j$ is not a loop, and $\gamma(F_j)=1$ if $F_j$ is a loop.  If $|F_j|$ is even, $x'_j$ and $x''_j$ are both in either $V_1$ or $V_2$, so that $\gamma(F_j)=1$.  If $|F_j|$ is odd, $x'_j$ and $x''_j$ are in different sets of the partition $(V_1,V_2)$, so that $a''(x'_j)a'(x''_j)=((1\pm i)/\sqrt2)^2$ and $\gamma(F_j)=\pm i$.  Choose now a pair $(G,\bar G)$ with odd $|G|=|\bar G|$ then $\gamma(G)+\gamma(\bar G)=0$ so that the terms in $p(z)$ corresponding to $F$ containing a connected component $G$ or $\bar G$ cancel.  This holds for all pairs $(G,\bar G)$ with odd $|G|=|\bar G|$ and therefore
$$	p(z)=\sum_{F\in U(E)}\prod_{F_j\subset F}^{\rm conn}\gamma(F_j)z^{|F_j|}
	=\sum_{F\in U_{\rm even}(E)}z^{|F|}=p_{\rm unbranched\,even}(z)   $$
\medskip
	With our choice of $a',a''$ we see that if $\alpha',\alpha''\in(-\pi/4,\pi/4)$ and
$$	{\rm Re}(z'_e+\epsilon)e^{-i\alpha'}>0\qquad,\qquad{\rm Re}(z''_e+\epsilon)e^{-i\alpha''}>0   $$
for all $e\in E$, we have $\prod_{x\in V}p_x^\epsilon(Z',Z'')\ne0$.  Therefore  by Asano contraction $P^\epsilon(Z)\ne0$ if
$$	-\pi/2<\alpha'+\alpha''<\pi/2\qquad{\rm and}\qquad(\forall e\in E)\,\,\,(z_e+\epsilon')e^{-i(\alpha'+\alpha'')}>0 $$
for some $\epsilon'>0$.  We take all $z_e$ equal to $z$ and use Hurwitz's theorem to let $\epsilon\to0$.  Since $\emptyset\in U_{\rm even}$, $p$ does not vanish identically and we obtain $p(z)\ne0$ if Re$(z)>0$, or by symmetry if Re$(z)\ne0$.\qed
\bigskip\noindent
{\bf 3. Oriented subgraphs of a non-oriented graph.}
\medskip
	Let $(V,E_0)$ be a non-oriented graph.  There are different ways to associate an oriented graph with $(V,E)$.  Here we define the oriented graph $(V,\tilde E_0)$ where each non-oriented edge $e\in E_0$ with endpoints $x_1,x_2\in V$ is replaced by two oriented edges $e',e''\in \tilde E_0$ such that $x'(e')=x_1,\,x''(e')=x_2$ and $x'(e'')=x_2,\,x''(e'')=x_1$.  We have thus $|\tilde E_0|=2|E_0|$.  The subgraphs $\tilde F$ of $(V,\tilde E_0)$, i.e., the subsets of $\tilde E_0$ may be called oriented subgraphs of $(V,E_0)$.
\medskip
	{\bf 3.1. Unbranched subgraphs of a non-oriented graph.}
\medskip
	From Proposition 2.1 we know that the polynomial counting oriented unbranched subgraphs of $(V,E_0)$, i.e.,
$$	p_{\rm oriented\,unbranched} (z)=\sum_{\tilde F\in U(\tilde E_0)}z^{|\tilde F|}
	=\sum_{\tilde F\in U(\tilde E_0)}\prod_{\tilde F_j\subset\tilde F}^{\rm conn}z^{|\tilde F_j|}   $$
has all its zeros real $<0$.  Note that without orientation
$$	p_{\rm unbranched} (z)=\sum_{F\in U(E_0)}z^{|F|}
	=\sum_{F\in U(E_0)}\prod_{F_j\subset F}^{\rm conn}z^{|F_j|}   $$
and it is known (see [5]) that this has all its zeros with real part $<0$.  [The set $U(E_0)$ of unbranched subgraphs of $E_0$ and the connected components of a non-oriented $F$ are defined in the obvious manner].  It is interesting to compare the oriented connected components $\tilde F_j$ with the non-oriented connected components $F_j$.  If $|F_j|=1$ then $F_j=e$ for some non-oriented $e\in E_0$ with endpoints $x_1,x_2$ and there are two oriented edges $e',e''\in\tilde E_0$ corresponding to $e$.  Also, corresponding to $F_j$ there are three connected components $\tilde F_{j\alpha}\subset\tilde E_0$, namely $\{e'\},\{e''\},\{e',e''\}$, and $|\tilde F_{j\alpha}|$ is 1 or 2.  If $|F_j|>1$, there correspond to $F_j$ two oriented $F_{j\alpha}$.  Assuming that $(V,E_0)$ has only simple edges between vertices we obtain thus for $\tilde E_0$ the polynomial
$$	p_{\rm oriented\,unbranched} (z)
	=\sum_{F\in U(E_0)}\prod_{F_j:|F_j|=1}^{\rm conn}(2z+z^2)\prod_{F_j:|F_j|>1}^{\rm conn}(2z^{|F_j|})   $$
\medskip
	{\bf 3.2. Even oriented unbranched subgraphs of a non-oriented graph.}
\medskip
	For a bipartite graph $E_0$ we obtain pairs $(G,\bar G)$ of subgraphs of $\tilde E_0$ as in Proposition 2.3 by orientation reversal so that
$$	p_{\rm oriented\,unbranched\,even} (z)=\sum_{\tilde F\in U_{\rm even}(\tilde E_0)}z^{|\tilde F|}   $$
has all its zeros purely imaginary by Proposition 2.3.
\medskip
	Let $(V,E_0)$ have only simple edges between vertices.  We define
$$	U'(E_0)=\{F:\hbox{ for all connected components $F_j$ of $F$ either $|F_j|=1$ or $|F_j|$ is even}\}   $$
Then we have
$$	p_{\rm oriented\,unbranched\,even} (z)=\sum_{F\in U'(E_0)}z^{2|\{|j:|F_j|=1\}|}.\prod_{j:|F_j|>1}(2z)^{|F_j|}   $$
for the unbranched even subgraph counting polynomial of $\tilde E_0$.
\bigskip\noindent
{\bf References.}
\medskip\noindent
[1] J.L. Lebowitz, B. Pittel, D. Ruelle, and E.R. Speer.  ``Central limit theorems, Lee-Yang zeros, and graph-counting polynomials.''  J. Combin. Theory Ser. A {\bf 142},147-183(2016).
\medskip\noindent
[2] J.L. Lebowitz, and D. Ruelle.  ``Phase transitions with four-spin interactions.''  Commun. Math. Phys. {\bf 311},755-768(2011).
\medskip\noindent
[3] J.L. Lebowitz, D. Ruelle, and E.R. Speer.  ``Location of the Lee-Yang zeros and absence of phase transitions in some Ising spin systems.''  J. Math. Phys. {\bf 53}095211(1-13)(2012).
\medskip\noindent
[4] D. Ruelle.  ``Zeros of graph-counting polynomials.''  Commun. Math. Phys. {\bf 200},43-56(1999).
\medskip\noindent
[5] D. Ruelle.  ``Counting unbranched subgraphs.''  J. Algebraic Combinatorics {\bf 9},157-160(1999).
\medskip\noindent
[6] D. Ruelle.  ``Characterization of Lee-Yang polynomials.''  Annals of Math. {\bf 171},589-603(2010).
\end